\documentclass[11pt]{amsart}
\usepackage{amssymb}
\usepackage{amsmath}
\usepackage{fancyhdr}
\usepackage[british]{babel}
\usepackage{geometry,mathtools}
\usepackage{enumitem}
\usepackage{algpseudocode}
\usepackage{dsfont}
\usepackage{centernot}
\usepackage{xstring}
\usepackage{colortbl}
\usepackage[section]{algorithm}
\usepackage{graphicx}
\usepackage[font={footnotesize}]{caption}
\usepackage[usenames,dvipsnames,table]{xcolor}
\usepackage{pstricks,tikz}
\usepackage[h]{esvect}
\usepackage[
		bookmarksopen=true,
		bookmarksopenlevel=1,
		colorlinks=true,
		linkcolor=darkblue,
        linktoc=page,
		citecolor=darkblue,
]{hyperref}
\usepackage{bbm}

\title[Hypergraph matchings]{Pseudorandom hypergraph matchings}
\author[S.~Ehard]{Stefan Ehard}
\address[S.~Ehard]{Institut f\"ur Optimierung und Operations Research, Universit\"at Ulm, Ulm,
Germany
}
\email{stefan.ehard@uni-ulm.de}

\author[S.~Glock]{Stefan Glock}
\author[F.~Joos]{Felix Joos}
\address[S.~Glock, F.~Joos]{School of Mathematics, University of Birmingham,
Edgbaston, Birmingham, B15 2TT, United Kingdom}
\email{[s.glock, f.joos]@bham.ac.uk}

\thanks{The research leading to these results was partially supported by the EPSRC, grant nos.~EP/N019504/1~(S.~Glock) and 
by the Deutsche Forschungsgemeinschaft (DFG, German Research Foundation) -- 339933727~(F.~Joos).
}

\newtheorem{theorem}[algorithm]{Theorem}

\newtheorem{lemma}[algorithm]{Lemma}

\theoremstyle{definition}

\newtheoremstyle{claimstyle}{5pt}{5pt}{\em}{5pt}{\em}{:}{5pt}{}
\theoremstyle{claimstyle}

\newtheoremstyle{stepstyle}{10pt}{5pt}{\em}{0pt}{\em}{:}{5pt}{}
\theoremstyle{stepstyle}
\newtheorem{step}{Step}

\numberwithin{equation}{section}

\definecolor{darkblue}{rgb}{0,0,0.5}

\def\noproof{{\unskip\nobreak\hfill\penalty50\hskip2em\hbox{}\nobreak\hfill%
       $\square$\parfillskip=0pt\finalhyphendemerits=0\par}\goodbreak}
\def\endproof{\noproof\bigskip}

\newdimen\margin
\def\textno#1&#2\par{
   \margin=\hsize
   \advance\margin by -4\parindent
          \setbox1=\hbox{\sl#1}
   \ifdim\wd1 < \margin
      $$\box1\eqno#2$$
   \else
      \bigbreak
      \hbox to \hsize{\indent$\vcenter{\advance\hsize by -3\parindent
      \it\noindent#1}\hfil#2$}
      \bigbreak
   \fi}

\def\proof{\removelastskip\penalty55\medskip\noindent\setcounter{claim}{0}\setcounter{step}{0}{\bf Proof. }} 

\def\lateproof#1{\removelastskip\penalty55\medskip\noindent\setcounter{claim}{0}\setcounter{step}{0}{\bf Proof of #1. }} 

\DeclareMathOperator{\dg}{deg}

\begin{document}

\newcommand{\new}[1]{\textcolor{red}{#1}}
\def\COMMENT#1{}
\def\TASK#1{}
\let\TASK=\footnote             

\newcommand{\todo}[1]{\begin{center}\textbf{to do:} #1 \end{center}}

\def\eps{{\varepsilon}}
\def\heps{{\hat{\varepsilon}}}

\newcommand{\ex}{\mathbb{E}}
\newcommand{\pr}{\mathbb{P}}
\newcommand{\cB}{\mathcal{B}}
\newcommand{\cA}{\mathcal{A}}
\newcommand{\cE}{\mathcal{E}}
\newcommand{\cS}{\mathcal{S}}
\newcommand{\cF}{\mathcal{F}}
\newcommand{\cG}{\mathcal{G}}
\newcommand{\bL}{\mathbb{L}}
\newcommand{\bF}{\mathbb{F}}
\newcommand{\bZ}{\mathbb{Z}}
\newcommand{\cH}{\mathcal{H}}
\newcommand{\cC}{\mathcal{C}}
\newcommand{\cM}{\mathcal{M}}
\newcommand{\bN}{\mathbb{N}}
\newcommand{\bR}{\mathbb{R}}
\def\O{\mathcal{O}}
\newcommand{\cP}{\mathcal{P}}
\newcommand{\cQ}{\mathcal{Q}}
\newcommand{\cR}{\mathcal{R}}
\newcommand{\cJ}{\mathcal{J}}
\newcommand{\cL}{\mathcal{L}}
\newcommand{\cK}{\mathcal{K}}
\newcommand{\cD}{\mathcal{D}}
\newcommand{\cI}{\mathcal{I}}
\newcommand{\cN}{\mathcal{N}}
\newcommand{\cV}{\mathcal{V}}
\newcommand{\cT}{\mathcal{T}}
\newcommand{\cU}{\mathcal{U}}
\newcommand{\cX}{\mathcal{X}}
\newcommand{\cZ}{\mathcal{Z}}
\newcommand{\cW}{\mathcal{W}}
\newcommand{\1}{{\bf 1}_{n\not\equiv \delta}}
\newcommand{\eul}{{\rm e}}
\newcommand{\Erd}{Erd\H{o}s}
\newcommand{\cupdot}{\mathbin{\mathaccent\cdot\cup}}
\newcommand{\whp}{whp }
\newcommand{\bX}{\mathcal{X}}
\newcommand{\bV}{\mathcal{V}}
\newcommand{\hbX}{\widehat{\mathcal{X}}}
\newcommand{\hbV}{\widehat{\mathcal{V}}}
\newcommand{\hX}{\widehat{X}}
\newcommand{\hV}{\widehat{V}}
\newcommand{\tX}{\widetilde{X}}
\newcommand{\tV}{\widetilde{V}}
\newcommand{\cbI}{\overline{\mathcal{I}^\alpha}}
\newcommand{\hAj}{\widehat{A}^\sigma_j}
\newcommand{\hVM}{V^M}
\newcommand{\supp}{{\rm supp}}

\newcommand{\doublesquig}{%
  \mathrel{%
    \vcenter{\offinterlineskip
      \ialign{##\cr$\rightsquigarrow$\cr\noalign{\kern-1.5pt}$\rightsquigarrow$\cr}%
    }%
  }%
}

\newcommand{\defn}{\emph}

\newcommand\restrict[1]{\raisebox{-.5ex}{$|$}_{#1}}

\newcommand{\cprob}[2]{\prob{#1 \;\middle|\; #2}}
\newcommand{\prob}[1]{\mathrm{\mathbb{P}}\left[#1\right]}
\newcommand{\expn}[1]{\mathrm{\mathbb{E}}\left[#1\right]}
\def\gnp{G_{n,p}}
\def\G{\mathcal{G}}
\def\lflr{\left\lfloor}
\def\rflr{\right\rfloor}
\def\lcl{\left\lceil}
\def\rcl{\right\rceil}

\newcommand{\qbinom}[2]{\binom{#1}{#2}_{\!q}}
\newcommand{\binomdim}[2]{\binom{#1}{#2}_{\!\dim}}

\newcommand{\grass}{\mathrm{Gr}}

\newcommand{\brackets}[1]{\left(#1\right)}
\def\sm{\setminus}
\newcommand{\Set}[1]{\{#1\}}
\newcommand{\set}[2]{\{#1\,:\;#2\}}
\newcommand{\krq}[2]{K^{(#1)}_{#2}}
\newcommand{\ind}[1]{$\mathbf{S}(#1)$}
\newcommand{\indcov}[1]{$(\#)_{#1}$}
\def\In{\subseteq}
\newcommand{\IND}{\mathbbm{1}}
\newcommand{\norm}[1]{\|#1\|}
\newcommand{\normv}[1]{\|#1\|_v}
\newcommand{\normV}[2]{\|#1\|_{#2}}

\newcommand{\Bij}{{\rm Bij}}
\newcommand{\inj}{{\rm inj}}
\newcommand{\mad}{{\rm mad}}

\begin{abstract} 
\noindent
A celebrated theorem of Pippenger states that any almost regular hypergraph with small codegrees has an almost perfect matching. We show that one can find such an almost perfect matching which is `pseudorandom', meaning that, for instance, the matching contains as many edges from a given set of edges as predicted by a heuristic argument.
\end{abstract}

\maketitle

\section{Introduction}
A hypergraph $\cH$ consists of a vertex set $V(\cH)$ and an edge set $E(\cH)\In 2^{V(\cH)}$. If all edges have size $r$, then $\cH$ is called \defn{$r$-uniform}, or simply an \defn{$r$-graph}. A \emph{matching} in $\cH$ is a collection of pairwise disjoint edges, and a \defn{cover of $\cH$} is a set of edges whose union contains all vertices. A matching is \defn{perfect} if it is also a cover.
These concepts are widely applicable, as `almost all combinatorial questions can be reformulated as either a matching or a covering problem of a hypergraph'~\cite{furedi:88}, and their study is thus of great relevance in combinatorics and beyond.

Results like Hall's theorem and Tutte's theorem that characterize when a graph has a perfect matching are central in graph theory.
However, for each $r\ge 3$, it is NP-complete to decide whether a given $r$-uniform hypergraph has a perfect matching. It is thus of great importance to find sufficient conditions that guarantee a perfect matching in an $r$-uniform hypergraph. This problem has received a lot of attention over the years. For instance, one line of research has focused on minimum degree conditions that guarantee a perfect matching (see e.g.~\cite{AGS:09,HPS:09,KM:15,RRS:09} and the survey~\cite{RR:10}). Another important direction has been to study
perfect matchings in random hypergraphs. The so-called Shamir's problem, to determine the threshold for which the (binomial) random $k$-graph has with high probability a perfect matching, was open for over 25~years resisting numerous efforts, until famously solved by Johansson, Kahn and Vu~\cite{JKV:08}. Moreover, Cooper, Frieze, Molloy and Reed~\cite{CFMR:96} determined when regular hypergraphs have with high probability a perfect matching.
It would be very interesting to obtain such results not only for random hypergraphs, but to find pseudorandomness conditions that (deterministically) guarantee a perfect matching. Apart from some partial results (e.g.~\cite{FK:12}), this seems wide open.

Many of the aforementioned results are proven by first obtaining an \emph{almost} perfect matching, and then using some clever ideas to complete it. It turns out that almost perfect matchings often exist under weaker conditions. For example, in the minimum degree setting, the threshold for finding an almost perfect matching is often smaller than that of finding a perfect matching. Also, there is a well-known theorem that yields almost perfect matchings under astonishingly mild pseudorandomness conditions. Mostly referred to as Pippenger's theorem, any almost regular hypergraph with small codegrees has an almost perfect matching. Both the result itself and also its proof method, the so-called `semi-random method' or `R\"odl nibble', have had a tremendous impact on Combinatorics.
We add to this body of research by showing the existence of `pseudorandom' matchings in this setting. We note that our result does not improve previous bounds on the \emph{size} of a matching that can be obtained. Rather, our focus is on the \emph{structure} of such a matching within the hypergraph it is contained in.

In Section~\ref{subsec:pippenger}, we revisit Pippenger's theorem. In Section~\ref{subsec:alon-yuster}, we discuss a theorem of Alon and Yuster, which can be viewed as an intermediate step. In Section~\ref{subsec:results}, we will motivate and state our main results.

\subsection{Pippenger's theorem} \label{subsec:pippenger}
Pippenger never published his theorem, and it was really the culmination of the efforts of various researchers in the 1980s. Most notably, in 1985, R\"odl~\cite{rodl:85} proved a long-standing conjecture of Erd\H{o}s and Hanani on approximate Steiner systems. A \defn{(partial) $(n,k,t)$-Steiner system} is a set $\cS$ of $k$-subsets of some $n$-set $V$ such that every $t$-subset of $V$ is contained in (at most) one $k$-set in~$\cS$. Steiner asked in 1853 for which parameters such systems exist, a question that has intrigued mathematicians for more than 150~years and was only answered recently by Keevash~\cite{keevash:14}. In 1963, Erd\H{o}s and Hanani asked whether one can, for fixed $k,t$, always find an `approximate Steiner system', that is, a partial $(n,k,t)$-Steiner system covering all but $o(n^t)$ of the $t$-sets, as $n\to \infty$. This was proved by R\"odl using the celebrated `nibble' method, with some ideas descending from~\cite{AKS:81,KPS:82}. Frankl and R\"odl~\cite{FR:85} observed that in fact a much more general theorem holds, which applies to almost regular hypergraphs with small codegrees. Pippenger's version stated below is a slightly stronger and cleaner version.
For a hypergraph $\cH$, we denote by $v(\cH)$ and $e(\cH)$ the number of vertices and edges of $\cH$, respectively, and we define for vertices $u,v\in V(\cH)$, the degree $\dg_{\cH}(v):=|\{e\in E(\cH)\colon v\in e \}|$ and codegree $\dg_{\cH}(uv):=|\{e\in E(\cH)\colon \{u,v \}\subseteq e \}|$.
Let
\begin{align*}
\Delta(\cH):=\max_{v\in V(\cH)}\dg_{\cH}(v), \quad
\delta(\cH):=\min_{v\in V(\cH)}\dg_{\cH}(v)
\text{~~~and~~~}
\Delta^c(\cH):=\max_{u\neq v\in V(\cH)} \dg_{\cH}(uv)
\end{align*}
denote the \defn{maximum degree}, \defn{minimum degree} and \defn{maximum codegree} of~$\cH$, respectively.

\begin{theorem}[Pippenger]\label{thm:pippenger}
For $r\in \bN$ and $\eps>0$, there exists $\mu>0$ such that any $r$-uniform hypergraph $\cH$ with $\delta(\cH)\ge (1-\mu)\Delta(\cH)$ and $\Delta^c(\cH)\le \mu \Delta(\cH)$ has a matching that covers all but at most an $\eps$-fraction of the vertices. 
\end{theorem}

To see why this generalizes R\"odl's result, fix $n,k,t$ and construct a hypergraph $\cH$ with vertex set $\binom{[n]}{t}$ where every $k$-set $X\In [n]$ induces the edge $\binom{X}{t}$. Note that perfect matchings in $\cH$ correspond to $(n,k,t)$-Steiner systems.
Clearly, $\cH$ is $\binom{k}{t}$-uniform. Moreover, every vertex has degree $\binom{n-t}{k-t}=\Theta(n^{k-t})$ and $\Delta^c(\cH)= \binom{n-t-1}{k-t-1}=o(n^{k-t})$. Thus, for sufficiently large~$n$, Pippenger's theorem implies the existence of a matching $\cM$ in $\cH$ that covers all but $o(n^t)$ of the vertices, which corresponds to a partial $(n,k,t)$-Steiner system which covers all but $o(n^t)$ of the $t$-sets. 
Frankl and R\"odl~\cite{FR:85} also applied (their version) of this theorem to obtain similar results for other combinatorial problems, for instance the existence of Steiner systems in vector spaces. Keevash~\cite{keevash:18c} raised the meta question of whether there exists a general theorem that provides sufficient conditions for a sparse `design-like' hypergraph to admit a perfect matching (for a notion of `design-like' that captures for example Steiner systems, but hopefully many more structures). Since such hypergraphs will likely be (almost) regular and have small codegree, the existence of an almost perfect matching follows from Pippenger's theorem, and a natural approach would be to use the absorbing method to complete such a matching to a perfect one. This of course can be extremely challenging since the relevant auxiliary hypergraphs are generally very sparse.

\subsection{The Alon--Yuster theorem} \label{subsec:alon-yuster}
In the case of Steiner systems, the absorbing method has been successfully applied to answer Steiner's question~\cite{GKLO:16,keevash:14}. Very roughly speaking, the idea of an absorbing approach is to set aside a `magic' absorbing structure, then to obtain an approximate Steiner system, and finally to employ the magic absorbing structure to clean up. One (minor, but still relevant) challenge is that the leftover of the approximate Steiner system must be `well-behaved'. More precisely, instead of the global condition that the number of uncovered $t$-sets is~$o(n^t)$, one needs the stronger local condition that for every fixed $(t-1)$-set, the number of uncovered $t$-sets containing this $(t-1)$-set is~$o(n)$.
Fortunately, Alon and Yuster~\cite[Theorem~1.2]{AY:05}, by building on a theorem of Pippenger and Spencer~\cite{PS:89}, provided a tool achieving this. They showed that any almost regular hypergraph with small codegrees contains a matching that is `well-behaved' in the sense that it not only covers all but a tiny proportion of the entire vertex set, but also has this property with respect to a specified collection of not too many not too small vertex subsets. In the above application to Steiner systems, for every $(t-1)$-set~$S$, consider the set $U_S\In V(\cH)$ of all $t$-sets containing~$S$. A matching in $\cH$ which covers almost all vertices of~$U_S$ then corresponds to a partial Steiner system which covers all but $o(n)$ of the $t$-sets containing~$S$, as desired.

\subsection{Pseudorandom matchings} \label{subsec:results}
The purpose of this paper is to provide a tool that is (qualitatively) a generalization of the Alon--Yuster theorem and gives much more control on the matching obtained. 
The need for such a tool arose in recent work of the authors on graph embeddings. In Section~\ref{sec:apps}, we will discuss further applications of our result in more detail.

To motivate this, suppose for simplicity that we are given a $D$-regular hypergraph and want to find an (almost) perfect matching~$\cM$. Moreover, we wish $\cM$ to be `pseudorandom', that is, to have certain properties that we expect from an idealized random matching. In a perfect matching, at a fixed vertex, exactly one edge needs to be included in the matching, and assuming that each edge is equally likely to be chosen, we may heuristically expect that every edge of $\cH$ is in a random perfect matching with probability~$1/D$. Thus, given a (large) set $E\In E(\cH)$ of edges, we expect $|E|/D$ matching edges in~$E$. 
More generally, given a set $X$, a \defn{weight function on~$X$} is a function $\omega\colon X\to \bR_{\ge 0}$.
For a subset $X'\In X$, we define $\omega(X'):=\sum_{x\in X'}\omega(x)$. If $\omega$ is a weight function on~$E(\cH)$, the above heuristic would imply that we expect from a `pseudorandom' matching $\cM$ that $\omega(\cM)\approx \omega(E(\cH))/D$. The following is a simplified version of our main theorem (Theorem~\ref{thm:hypermatching2}) which asserts that a hypergraph with small codegrees has a matching that is pseudorandom in the above sense.

\begin{theorem}\label{thm:hypermatching simple}
Suppose $\delta\in(0,1)$ and $r\in \bN$ with $r\ge 2$, and let $\eps:=\delta/50r^2$. Then there exists $\Delta_0$ such that for all $\Delta\ge \Delta_0$, the following holds:
Let $\cH$ be an $r$-uniform hypergraph with $\Delta(\cH)\leq \Delta$ and $\Delta^c(\cH)\le \Delta^{1-\delta}$ as well as $e(\cH)\leq \exp(\Delta^{\eps^2})$. 
Suppose that $\cW$ is a set of at most $\exp(\Delta^{\eps^2})$ weight functions on~$E(\cH)$.
Then, there exists a matching $\cM$ in~$\cH$ such that $\omega(\cM)=(1\pm \Delta^{-\eps}) \omega(E(\cH))/\Delta$ for all $\omega \in \cW$ with $\omega(E(\cH))\ge \max_{e\in E(\cH)}\omega(e)\Delta^{1+\delta}$.
\end{theorem}

Let us discuss a few aspects of this theorem. First, note that we do not require $\cH$ to be almost regular. The theorem can be applied with any (sufficiently large) $\Delta$, and we will discuss the usefulness of this and the fact that $v(\cH)$ plays no role in the parametrization of the theorem, in more detail in Section~\ref{sec:apps}. If $\cH$ is almost regular, an almost perfect matching can be obtained by applying the theorem with $\Delta=\Delta(\cH)$ to the weight function $\omega\equiv 1$. This yields that $|\cM| \ge (1-o(1))\frac{e(\cH)}{\Delta(\cH)} \ge (1-o(1))v(\cH)/r$, where the last inequality uses that $re(\cH)=\sum_{x\in V(\cH)}\dg_{\cH}(x)=(1\pm o(1))v(\cH)\Delta(\cH)$.\COMMENT{(StefanE): Replace $\Delta(\cH)$ with $\Delta$? As we apply Theorem~\ref{thm:hypermatching simple} with $\Delta\geq \Delta(\cH)$ satisfying $\Delta^c(\cH)\le \Delta^{1-\delta}$.}

We remark that, while Pippenger's theorem only needs $\Delta^c(\cH)= o(\Delta)$, we need a stronger condition to apply concentration inequalities. For the same reason, we also need that $\omega(E(\cH))$ is not too small (relative to the maximum possible weight). As a result, our theorem also allows stronger conclusions in that the error term $\Delta^{-\eps}$ decays polynomially with~$\Delta$. 

Note that Theorem~\ref{thm:hypermatching simple} is (qualitatively) more general than the Alon--Yuster theorem. Indeed, suppose $\cH$ is an almost regular hypergraph and we are given a collection $\cV$ of subsets $U\In V(\cH)$ and want to ensure that $\cM$ covers each $U\in \cV$ almost completely. For each target subset $U\in \cV$, we can define a weight function $\omega_U$ by setting $\omega_U(e):=|e\cap U|$. Note that $\omega_U(E(\cH)) = \sum_{x\in U}\dg_\cH(x) = (1\pm o(1))|U|\Delta(\cH)$. Thus, if $\omega_U(\cM)=(1\pm o(1)) \omega_U(E(\cH))/\Delta(\cH)$, we deduce that $|U\cap V(\cM)|=\omega_U(\cM)=(1\pm o(1)) \omega_U(E(\cH))/\Delta(\cH) \ge (1-o(1))|U|$, implying that almost all vertices of $U$ are covered by~$\cM$.

\bigskip
In fact, we prove a more general theorem which not only allows weight functions on edges, but on tuples of edges. This allows, for instance, to specify a set of pairs of edges, and control how many pairs will be contained in the matching. In Section~\ref{subsec:pseudo Steiner}, we will use this to count subconfigurations of Steiner systems.

Given a set $X$ and an integer $\ell\in \bN$, 
an \defn{$\ell$-tuple weight function on~$X$} is a function $\omega\colon \binom{X}{\ell}\to \bR_{\ge 0}$, that is, a weight function on~$\binom{X}{\ell}$.
For a subset $X'\In X$, we then define $\omega(X'):=\sum_{S\in \binom{X'}{\ell}}\omega(S)$. Moreover, if $\cX\In \binom{X}{\ell}$, we write $\omega(\cX)$ for $\sum_{S\in \cX}\omega(S)$ as for usual weight functions. For $k\in[\ell]_0$ and a tuple $T\in \binom{X}{k}$, define 
\begin{align}
	\omega(T):=\sum_{S\supseteq T}\omega(S),\mbox{ and let }\normV{\omega}{k}:=\max_{T\in \binom{X}{k}}\omega(T). \label{def shadow weight}
\end{align}
Suppose $\cH$ is an $r$-uniform hypergraph and $\omega$ is an $\ell$-tuple weight function on~$E(\cH)$. Clearly, if $\cM$ is a matching, then a tuple of edges which do not form a matching will never contribute to~$\omega(\cM)$. We thus say that $\omega$ is \defn{clean} if $\omega(\cE)=0$ whenever $\cE\in\binom{E(\cH)}{\ell}$ is not a matching.

The following is our main result, which readily implies Theorem~\ref{thm:hypermatching simple}.\COMMENT{$\ell=L=1$. Since $\normV{\omega}{1}=\max_{e\in E(\cH)}\omega(e)$. Also, for $\ell=1$, everything is clean.}

\begin{theorem}\label{thm:hypermatching2}
Suppose $\delta\in(0,1)$ and $r\in \bN$ with $r\ge 2$, and let $\eps:=\delta/50L^2r^2$. Then there exists $\Delta_0$ such that for all $\Delta\ge \Delta_0$, the following holds:
Let $\cH$ be an $r$-uniform hypergraph with $\Delta(\cH)\leq \Delta$ and $\Delta^c(\cH)\le \Delta^{1-\delta}$ as well as $e(\cH)\leq \exp(\Delta^{\eps^2})$. 
Suppose that for each $\ell\in[L]$, we are given a set $\cW_\ell$ of clean $\ell$-tuple weight functions on $E(\cH)$ of size at most $\exp(\Delta^{\eps^2})$, such that $\omega(E(\cH))\ge \normV{\omega}{k}\Delta^{k+\delta}$ for all $\omega\in \cW_\ell$ and $k\in [\ell]$.

Then, there exists a matching $\cM$ in~$\cH$ such that $\omega(\cM)=(1\pm \Delta^{-\eps}) \omega(E(\cH))/\Delta^\ell$ for all $\ell\in [L]$ and $\omega \in \cW_\ell$.
\end{theorem}

We will prove Theorem~\ref{thm:hypermatching2} in Section~\ref{sec:proof}, after stating some preliminary results in the next section. In Section~\ref{sec:apps}, we will discuss applications of our main result.

\section{Preliminaries}

Our main tool is the next theorem of Molloy and Reed on the chromatic index of a hypergraph with small codegrees, improving on earlier work of Pippenger and Spencer as well as Kahn.
Pippenger and Spencer~\cite{PS:89} strengthened Theorem~\ref{thm:pippenger} by showing that under the same assumptions, one can even obtain an almost optimal edge-colouring of $\cH$, using $(1+o(1))\Delta$ colours. (The existence of an almost perfect matching follows then by averaging over the colour classes.)
Kahn~\cite{kahn:96} generalized this to list colourings, and Molloy and Reed improved the $o(1)$-term. For simplicity, we only state their result for normal colourings.

\begin{theorem}[Molloy and Reed~{\cite[Theorem~2]{MR:00}}]\label{thm:MR}
Let $1/\Delta\ll \delta, 1/r$.
Suppose $\mathcal{H}$ is an $r$-uniform hypergraph satisfying $\Delta^c(\mathcal{H})\leq \Delta^{\delta}$ and $\Delta(\mathcal{H})\le \Delta$.
Then, the edge set $E(\cH)$ can be decomposed into $\Delta+ \Delta^{1-\frac{1-\delta}{r}}\log^5 \Delta$ edge-disjoint matchings.\COMMENT{Original: $B=\Delta^c$, then need $\Delta+c_r \Delta^{1-1/r}B^{1/r}(\log \Delta /B)^4$. Can assume $B\ge 1$ and thus upper bound second term by $c_r \Delta^{1-1/r}B^{1/r}\log^4 \Delta$. Then use upper bound for $B$ to get $\Delta^{1-1/r}B^{1/r}\le \Delta^{1-\frac{1-\delta}{r}}$ and incorporate constant in $\log$-factor.}
\end{theorem}

Note here that $\cH$ is not required to be almost regular. In fact, this assumption can also be omitted from the Pippenger--Spencer theorem since any given $r$-uniform hypergraph $\cH$ can be embedded into a $\Delta(\cH)$-regular hypergraph $\cH'$ with $\Delta^c(\cH')=\Delta^c(\cH)$,\COMMENT{To obtain $\cH'$, proceed as follows: Take $r$ copies of $\cH$. For all vertices in $\cH$ of degree less than $\Delta$, add an edge on the $r$ copies of that vertex. Clearly, the codegree does not increase. Then repeat. The minimum degree increases in each step, thus after at most $\Delta$ steps, we have a $\Delta$-regular hypergraph.}
and any colouring of $\cH'$ induces a colouring of $\cH$ with the same number of colours. 

\bigskip
We also make use of several probabilistic tools to establish concentration of a random variable~$X$. If $X$ is the sum of independent Bernoulli variables, we use the following well-known Chernoff-type bound.

\begin{theorem}[Chernoff's bound]\label{thm:chernoff}
Suppose $X_1,\dots,X_m$ are independent random variables taking values in~$\Set{0,1}$. Let $X:=\sum_{i=1}^m X_i$. Then, for all $\lambda> 0$, $$\prob{|X-\expn{X}|\ge \lambda} \le 2 \exp\left({-\frac{\lambda^2}{2(\expn{X} + \lambda/3)}}\right).$$
\end{theorem}
\COMMENT{For $0\le \eps\le 3/2$, this implies (with $\lambda=\eps\expn{X}$) that $\prob{X\neq (1\pm \eps)\expn{X}} \le 2 \exp\left({-\eps^2\expn{X}/3}\right).$ If $\lambda\ge 7\expn{X}$, we obtain $\prob{X \ge \lambda} \le 2 \exp(-\lambda).$ Moreover, if $\expn{X}\le 1.1\lambda$, then $\prob{|X-\expn{X}|\ge \lambda} \le 2 \exp\left({-\frac{\lambda}{3}}\right)$}

Similarly, if $X$ is a function of several independent Bernoulli variables and does not depend too much on any of the variables, we use the following `bounded differences inequality'.

\begin{theorem}[McDiarmid's inequality, see~\cite{mcdiarmid:89}\COMMENT{Lemma~1.2}] \label{thm:McDiarmid}
Suppose $X_1,\dots,X_m$ are independent Bernoulli random variables and suppose $b_1,\dots,b_m\in [0,B]$.
Suppose $X$ is a real-valued random variable determined by $X_1,\dots,X_m$ such that changing the outcome of $X_i$ changes $X$ by at most $b_i$ for all $i\in [m]$.
Then, for all $\lambda>0$, we have $$\prob{|X-\expn{X}|\ge \lambda} \le 2 \exp\left({-\frac{2\lambda^2}{B\sum_{i=1}^m b_i}}\right).$$
\end{theorem}

In one of our proofs we consider exposure martingales;
that is, suppose we have a random variable $X$ that is determined by independent random variables $Y_1,\ldots,Y_n$
and we define $X_t:=\expn{X\mid Y_1,\ldots,Y_t}$.
Then it is well-known that $(X_t)_{t\geq 0}$ is a martingale, the so-called \defn{exposure martingale} for~$X$. Note that $X_0=\expn{X}$ and $X_n=X$. 
Now, Freedman's martingale concentration inequality can be used to obtain concentration of $X$ around its mean.

\begin{lemma}[Freedman's inequality~\cite{freedman:75}]\label{freedman}
Let $(\Omega, \cF, \mathbb{P})$ be a probability space and let $(\cF_t)_{t\geq 0}$ be a filtration of $\cF$.
Let $(X_t)_{t\geq 0}$ be a martingale adapted to $(\cF_t)_{t\geq 0}$.
Suppose $\sum_{t\ge 0} \expn{|X_{t+1}-X_{t}| \mid \cF_t} \le \sigma$ and that $|X_{t+1}-X_t|\leq C$ for all~$t$.
Then, for any $\lambda>0$,
\begin{align*}
\prob{|X_t-X_0| \geq \lambda \text{ for some }t}\leq 2\exp\left({-\frac{\lambda^2}{2C(\lambda + \sigma)}}\right).
\end{align*} 
\end{lemma}

For $a,b,c\in \mathbb{R}$,
we write $a=b\pm c$ whenever $a\in [b-c,b+c]$.
For $a,b,c\in (0,1]$,
we sometimes write $a\ll b \ll c$ in our statements meaning that there are increasing functions $f,g:(0,1]\to (0,1]$
such that whenever $a\leq f(b)$ and $b \leq g(c)$,
then the subsequent result holds.
We assume that large numbers are integers if that does not affect the argument.

\section{Proof} \label{sec:proof}

We first sketch our proof. For simplicity, we consider the setting of Theorem~\ref{thm:hypermatching simple}.
We split $\cH$ randomly into $p$ vertex-disjoint induced subgraphs $\cH_1,\dots,\cH_p$ and let $\cH'$ be the union of those. With high probability, $\Delta(\cH_i)\approx \Delta(\cH)p^{-(r-1)}$ for each $i$, and for a given weight function $\omega$, we have $\omega(E(\cH'))\approx \omega(E(\cH))p^{-(r-1)}$. After fixing such a partition, we utilize the theorem of Molloy and Reed to find, for each $i\in[p]$, a partition of $E(\cH_i)$ into $M\approx \Delta(\cH)p^{-(r-1)}$ matchings. Finally, we select a matching from each partition uniformly at random, and let $\cM$ be the union of these matchings. Clearly, every edge in $\cH'$ is contained in $\cM$ with probability~$M^{-1}$, so $\expn{\omega(\cM)}=\omega(E(\cH'))M^{-1}\approx \omega(E(\cH))/\Delta(\cH)$. Moreover, the individual effect of the matching chosen in $\cH_i$ is relatively small, so we can use McDiarmid's inequality to establish concentration.
This approach is the same as taken by Alon and Yuster. 
One important new ingredient in our proof is that we partition each $\cH_i$ further into edge-disjoint subgraphs $\cH_{i,1},\dots,\cH_{i,q}$ such that $\omega(E(\cH_{i,j}))$ is of magnitude $\omega(E(\cH_{i}))/q$, and then apply Theorem~\ref{thm:MR} to each~$\cH_{i,j}$. This gives, as above, a partition of $\cH_i$ into matchings, from which we still choose one uniformly at random. However, the individual effect of each matching chosen has now been drastically reduced, which allows us to apply McDiarmid's inequality with the desired parameters.

\medskip

In the setting of Theorem~\ref{thm:hypermatching simple}, the partition of each $\cH_i$ into edge-disjoint subgraphs $\cH_{i,1},\dots,\cH_{i,q}$ could be done easily with a generalized Chernoff bound. However, in the setting of Theorem~\ref{thm:hypermatching2}, we are not aware of a conventional concentration inequality that suits our needs for this step (in particular, since $q$ is rather large). Thus, we first prove a tool that will achieve this for us.
Roughly speaking, what we require is the following: Let $\cH$ be a `directed' $\ell$-graph on~$V$, that is, a collection of ordered $\ell$-subsets of~$V$. Let $f\colon V\to [q]$ be obtained by choosing $f(v)\in [q]$ uniformly at random for each vertex $v$ independently. For each directed edge $e=(v_1,\dots,v_\ell)$, let $f(e):=(f(v_1),\dots,f(v_\ell))$. 
For a fixed `pattern' $\alpha \in [q]^\ell$, let $X_\alpha$ denote the number of $e\in E(\cH)$ with $f(e)=\alpha$. 
Clearly, for each edge~$e$, we have that $\prob{f(e)=\alpha}=q^{-\ell}$, thus, $\expn{X_\alpha}=q^{-\ell}e(\cH)$. We would like to know that $X_\alpha$ is concentrated around its mean, even when $q$ is quite large.

For simplicity, we will actually only consider the case when $\cH$ is an $\ell$-graph, the vertex set $V$ is ordered, and each edge of $\cH$ obtains its direction from the ordering of~$V$.
Thus, our setup is as follows. Let $(V,<)$ be an ordered set. Let $f\colon V\to [q]$ be obtained by choosing $f(v)\in [q]$ uniformly at random for each $v\in V$ independently. For each $\ell$-set $e=\Set{v_1,\dots,v_\ell}$ with $v_1<\dots<v_\ell$, let $f(e):=(f(v_1),\dots,f(v_\ell))$. For a fixed `pattern' $\alpha \in [q]^\ell$, let $E_\alpha=E_\alpha(f)$ denote the (random) set of all $e\in \binom{V}{\ell}$ with $f(e)=\alpha$. Given an $\ell$-tuple weight function $\omega$ on~$V$, the following theorem shows that the random variable $\omega(E_{\alpha})$ is concentrated around its mean. 

\begin{theorem}\label{thm: edge slicing}
Suppose $(V,<)$, $f$, $\alpha$, $\omega$ are as above. Suppose that $g\ge 24\ell^3 (\ell+1+\log |V|)$. Define $M:= q^{-\ell}\max_{k\in[\ell]}\{\normV{\omega}{k}q^{k}g^{k-1}\}$.
Then for any $\lambda> 0$, we have
\begin{align*}
\prob{|\omega(E_\alpha)-  \expn{\omega(E_\alpha)}|\ge \lambda }
\leq 2^{\ell} \exp\left(- \frac{\lambda^2}{12\ell^2M(\lambda+\expn{\omega(E_\alpha)})}\right) + \exp\left(-\frac{g}{24\ell^2}\right).
\end{align*}
\end{theorem}

\proof
Let $n:=|V|$ and let $v_1<\dots<v_n$ be the ordered elements of~$V$\COMMENT{could also just take $[n]$ and then have $j_i$ instead of $v_{j_i}$} and write $\alpha=(\alpha_1,\dots,\alpha_\ell)$. For $t\in [n]_0$, let $$X_t:=\expn{\omega(E_\alpha)\mid f(v_1),\ldots, f(v_t)}$$ (and $X_t:=X_n$ for $t\geq n$).
Hence $X=(X_t)_{t\geq 0}$ is the so-called exposure martingale for~$\omega(E_\alpha)$, where the labels $f(v_i)$ are revealed one by one. In particular, $X_0=\expn{\omega(E_\alpha)}$ and $X_n=\omega(E_\alpha)$.

For $k\in [\ell]$ and a $k$-tuple weight function $\omega'$ on $V$, let $$M_k(\omega'):= q^{-k} \max_{i\in [k]}\{ \normV{\omega'}{i}q^{i}g^{i-1}\}.$$ Note that we have
\begin{align}
	M_k(\omega') q^k &\le M_\ell(\omega') q^\ell. \label{norm monotonicity}
\end{align}
Let $M_k:=M_k(\omega)$ and note that $M=M_\ell$.

We prove the theorem by induction on~$\ell$ (with $(V,<)$ and $g$ being fixed). Thus, assume first that $\ell=1$. (This case is also contained in the inductive step below with no inductive hypothesis being needed, but the short proof here may serve as a warm up.) Observe that $X_t(f)-X_{t-1}(f)=\omega(\Set{v_t})(\IND_{f(v_t)=\alpha_1}-1/q)$ for $t\in[n]$. Hence, we can directly apply Freedman's inequality\COMMENT{$|X_t-X_{t-1}|\le \omega(\Set{v_t}) \le \normV{\omega}{1}$ and $\expn{|X_t-X_{t-1}| \mid f(v_1),\ldots,f(v_{t-1})}= 2\omega(\Set{v_t})(1-1/q)q\le 2\omega(\Set{v_t})/q$} to obtain (observe that $M_1=\normV{\omega}{1}$)
\begin{align*}
\prob{|\omega(E_\alpha)-  \expn{\omega(E_\alpha)}|\ge \lambda }
&\leq 2\exp\left(- \frac{\lambda^2}{2\normV{\omega}{1}(\lambda + \sum_{t\in [n]}2\omega(\Set{v_t})/q )}\right)\\
&\le 2\exp\left(- \frac{\lambda^2}{4 M_1(\lambda + \expn{\omega(E_\alpha)}) }\right),
\end{align*}
as desired.

Suppose now that $\ell\geq 2$.
In order to apply induction, we need to introduce some more notation.
For $t\in [n]$ and $k\in [\ell-1]_0$, let $\omega^{t,k}\colon \binom{V}{k}\to [0,\infty)$ be defined as (where $j_1<\ldots <j_k$)
$$\omega^{t,k}(\{v_{j_1},\ldots,v_{j_{k}}\})
:=\sum_{\substack{j_{k+1} < \ldots< j_\ell\\ j_k < j_{k+1}=t}}
\omega(\{v_{j_1},\ldots,v_{j_\ell}\}).$$
Moreover, let $\omega^{\leq t,k}\colon \binom{V}{k}\to [0,\infty)$ be defined by $\omega^{\leq t,k}(S):=\sum_{s\leq t}\omega^{s,k}(S)$ for all $S\in \binom{V}{k}$.
Note that
\begin{align}
	\mbox{$\omega^{\leq n,k}(V) = \omega(V)$ and $\normV{\omega^{\leq n,k}}{i} \le \normV{\omega}{i}$ for all $i\in[k]$.}\label{norm inherit}
\end{align}
\COMMENT{Observe that every $\ell$-edge $(v_1,\ldots,v_\ell)$ transfers its weight to $(v_1,\ldots,v_k)$ through the weight function $\omega^{t,k}$ where $t=k+1$.
This also explains the second part easily.
}

For $k\in[\ell-1]_0$, let $\alpha[k]:=(\alpha_1,\dots,\alpha_k)$, and define $E_{\alpha[k]}=E_{\alpha[k]}(f)$ as the random set of all $k$-sets $\Set{v_{j_1},\dots,v_{j_k}}$ for which $f(v_{j_i})=\alpha_i$ for all $i\in[k]$, where $j_1<\dots<j_k$. For clarity, we briefly discuss the case $k=0$, when $\omega^{t,0}$ is the function that maps $\emptyset$ to $\sum_{t < j_2 <\ldots< j_\ell} \omega(\{v_{t},v_{j_2}\ldots,v_{j_\ell}\})$. In particular, we have for all $t\in[n]$ that
\begin{align}
\omega^{t,0}(\emptyset) &\le \omega(\Set{v_t}) \le \normV{\omega}{1}=M_1;\label{induction hypothesis null 1} \\
\omega^{\le t,0}(\emptyset) &\le \omega(V).\label{induction hypothesis null 2}
\end{align}
Note also that $E_{\alpha[0]}=\Set{\emptyset}$.

The purpose of these definitions lies in the following formula for the one-step change of the process~$X$: for $t\in[n]$, we have
\begin{align*}
X_t(f)-X_{t-1}(f) &= \sum_{k\in[\ell-1]_0} \omega^{t,k}(E_{\alpha[k]}(f))\cdot (\IND_{f(v_t)=\alpha_{k+1}}-1/q)\cdot q^{-(\ell-(k+1))}. 
\end{align*}
\COMMENT{Split the sum over all $\ell$ tuples. Those which do not contain $v_t$ cancel out. For those which contain $v_t$, split up into how many vertices $k$ they contain left of~$v_t$. For those, there are still $\ell-(k-1)$ vertices right from $v_t$ each attaining the required label with probability $1/q$. The only difference between $X_t$ and $X_{t-1}$ is whether $v_t$ itself attains the correct label. In $X_{t-1}$, $f(v_t)$ is not yet fixed, so we get $1/q$. In $X_t$, $f(v_t)$ is already fixed, so we get $\IND_{f(v_t)=\alpha_{k+1}}$}%
Clearly, $|\IND_{f(v_t)=\alpha_{k+1}}-1/q|\le 1$ and $\expn{|\IND_{f(v_t)=\alpha_{k+1}}-1/q|}=2(1-1/q)/q\le 2/q$. Hence, for the absolute change and expected absolute change of the process~$X$ in one step we obtain the following bounds:
\begin{align}
|X_t-X_{t-1}| &\leq \sum_{k\in[\ell-1]_0}\omega^{t,k}(E_{\alpha[k]})\cdot q^{k+1 -\ell}; \label{original absolute change} \\
\expn{|X_t-X_{t-1}| \mid f(v_1),\ldots,f(v_{t-1})} &\le \sum_{k\in[\ell-1]_0} 2\omega^{t,k}(E_{\alpha[k]}) \cdot q^{k-\ell}.\label{original expected change}
\end{align}
Note that $\omega^{t,k}(E_{\alpha[k]})$ is itself a random variable, when $k>0$. Unfortunately, its deterministic upper bound is not good enough to apply Freedman's inequality directly to the martingale $(X_t)_{t\ge 0}$. We apply a common trick by defining a stopped process $Y=(Y_t)_{t\ge 0}$ which is equal to $X$ as long as the random variables $\omega^{t,k}(E_{\alpha[k]})$ behave nicely, and then `freezes'. We can then apply Freedman's inequality to~$Y$. Finally, we need to show that the process is unlikely to freeze, implying that the concentration result for $Y$ transfers to~$X$. For this, we employ the statement inductively with $\omega^{t,k},\omega^{\le n,k},\alpha[k]$.

We define two types of stopping times for~$X$.
For $k\in [\ell-1]$, let 
\begin{align}
\tau'_k:= \min_{t\in [n-1]}\{\omega^{\le t+1,k}(E_{\alpha[k]})\geq \omega(V)q^{-k}+\lambda q^{\ell-k}\}\wedge n. \label{stopping time 1}
\end{align}
Moreover, for $k\in [\ell-1]$ and $t\in[n-1]$, define
\begin{align}\label{stopping time 2}
\tau^t_k:=
\begin{cases}
t &\text{if $\omega^{t+1,k}(E_{\alpha[k]})\geq 2M_{k+1}$,} 
\\
n &\text{otherwise.}
\end{cases}
\end{align}
Let $\tau :=\min_{t\in [n],k\in [\ell-1]}\{\tau'_k,\tau_k^t\}$.
Note that $\omega^{t+1,k}(E_{\alpha[k]})$ is fully determined by $f(v_1),\dots,f(v_{t})$, since $\omega^{t+1,k}(S)=0$ whenever $S$ contains a vertex $v_j$ with $j\ge t+1$.  Thus, $\tau$ is indeed a stopping time for~$X$. We define $Y=(Y_t)_{t\geq 0}$ by $Y_t:=X_{t\wedge \tau}$, and let $\Delta Y_t:=Y_t-Y_{t-1}$.\COMMENT{We want that $X_n=Y_n$, so $\tau \ge n$ is enough. Works nicer with indices as well. Otherwise, we would often need to write $\tau\wedge n$} By the optional stopping theorem $Y$ is also a martingale, and thus we can apply Freedman's inequality. To this end, we next bound the absolute and expected one step change for~$Y$.

We claim that $|\Delta Y_t|\leq 2\ell M_\ell$ for all~$t$.
Indeed, if $t\ge \tau+1$, then trivially $|\Delta Y_t|=0$ and whenever $t\le \tau$, then
\begin{align*}
|\Delta Y_t|
&\overset{\eqref{original absolute change}}{\leq} \sum_{k\in[\ell-1]_0}\omega^{t,k}(E_{\alpha[k]})\cdot q^{k+1 -\ell} \overset{\eqref{induction hypothesis null 1},\eqref{stopping time 2}}{\leq} \sum_{k\in[\ell-1]_0}2M_{k+1} \cdot q^{k+1 -\ell} \overset{\eqref{norm monotonicity}}{\leq} 2\ell M_\ell.
\end{align*}

Similarly,
\begin{eqnarray*}
\sum_{t\ge 1}\expn{|\Delta Y_t| \mid f(v_1),\ldots,f(v_{t-1})}
&\overset{\eqref{original expected change}}{\leq}& \sum_{t\in [\tau]} \sum_{k\in[\ell-1]_0} 2\omega^{t,k}(E_{\alpha[k]}) \cdot q^{k-\ell} \\
&=& \sum_{k\in[\ell-1]_0} 2\omega^{\le \tau,k}(E_{\alpha[k]})\cdot q^{k-\ell} \\
&\overset{\eqref{induction hypothesis null 2},\eqref{stopping time 1}}{\leq}& \sum_{k\in[\ell-1]_0} 2(\omega(V)q^{-k}+\lambda q^{\ell-k}) \cdot q^{k-\ell} \\
&=&2\ell (\omega(V) q^{-\ell}+\lambda).
\end{eqnarray*}

Thus, we can apply Freedman's inequality to obtain
\begin{align*}
\prob{|Y_n-Y_0|\geq \lambda} &\leq 2\exp\left(- \frac{\lambda^2}{4\ell M_\ell (\lambda + 2\ell (\omega(V) q^{-\ell}+\lambda))}\right) \\
&\leq  2\exp\left(- \frac{\lambda^2}{12\ell^2 M_\ell(\lambda + \expn{\omega(E_\alpha)})}\right) .
\end{align*}

It remains to show that $Y_n=X_n$ with high probability.
We first consider the stopping times~$\tau'_k$. Fix $k\in [\ell-1]$ and note that $\expn{\omega^{\le n,k}(E_{\alpha[k]})} = \omega^{\le n,k}(V)/q^k = \omega(V)/q^k$ by~\eqref{norm inherit}. We apply the induction hypothesis to $\omega^{\leq n,k}$, with $\lambda q^{\ell-k}$ and $k$ playing the roles of $\lambda$ and $\ell$, and obtain
\begin{align*}
\prob{\tau'_k<n} 
&\leq \prob{\omega^{\le n,k}(E_{\alpha[k]}) \ge \expn{\omega^{\le n,k}(E_{\alpha[k]})} + \lambda q^{\ell-k} }\\
&\leq 2^{k}\exp\left( - \frac{\lambda^2 q^{2(\ell-k)}}{12k^2M_{k}(\omega^{\leq n,k}) (\lambda q^{\ell-k} + \expn{\omega^{\le n,k}(E_{\alpha[k]})})} \right) + \exp\left(-\frac{g}{24k^2}\right)\\
& \leq  2^{k}\exp\left(- \frac{\lambda^2}{12k^2 M_\ell(\lambda + \expn{\omega(E_\alpha)})}\right) + \exp\left(-\frac{g}{24k^2}\right),
\end{align*}
where we have used that $\expn{\omega^{\le n,k}(E_{\alpha[k]})}= q^{\ell-k}\expn{\omega(E_\alpha)}$ and $M_{k}(\omega^{\leq n,k})\le M_k(\omega)\le q^{\ell-k}M_\ell$ by~\eqref{norm inherit} and~\eqref{norm monotonicity}.

Next, we consider the stopping times $\tau_k^t$. Let $k\in [\ell-1]$ and $t\in[n-1]$.
Observe that $\normV{\omega^{t,k}}{i}\leq \normV{\omega}{i+1}$ for all $i\in [k]$.\COMMENT{This holds as $v_t$ is `contained' in every edge.}
Hence 
\begin{align*}
\frac{M_{k+1}(\omega)}{M_{k}(\omega^{t,k})}
&= \frac{q^{-k-1}\max_{i\in [k+1]}\{ \normV{\omega}{i}q^{i}g^{i-1}\}}
{q^{-k}\max_{i\in [k]}\{ \normV{\omega^{t,k}}{i}q^{i}g^{i-1}\}} 
\ge \frac{g\max_{i\in [k+1]}\{ \normV{\omega}{i}q^{i}g^{i-1}\}}{\max_{i\in [k+1]\setminus\{1\}}\{ \normV{\omega}{i}q^{i}g^{i-1}\}} \ge g.
\end{align*}
Note that $\expn{\omega^{t,k}(E_{\alpha[k]})}= q^{-k}\omega^{t,k}(V) \leq q^{-k}\normV{\omega}{1} \leq M_{k+1}$.\COMMENT{e.g. by~\eqref{norm monotonicity}}
Thus, using induction for $\omega^{t,k}$ with $M_{k+1}$ and $k$ playing the roles of $\lambda$ and $\ell$, we deduce that
\begin{align*}
\prob{\tau_k^t<n} 
&\leq \prob{\omega^{t,k}(E_{\alpha[k]})\geq 2M_{k+1}} \leq 2^{k}\exp\left( - \frac{M_{k+1}}{24 k^2 M_{k}(\omega^{t,k})} \right)+\exp\left(-\frac{g}{24k^2}\right) \\&\leq (2^{k}+1)\exp\left(-\frac{g}{24k^2}\right).
\end{align*}
A union bound now implies that
\begin{align*}
	\prob{\tau<n} &\le \sum_{k=1}^{\ell-1} \left(2^{k}\exp\left(- \frac{\lambda^2}{12k^2 M_\ell(\lambda + \expn{\omega(E_\alpha)})}\right) + (1+n(2^{k}+1)) \exp\left(-\frac{g}{24k^2}\right) \right) \\
	&\le (2^{\ell}-2) \exp\left(- \frac{\lambda^2}{12\ell^2 M_\ell(\lambda + \expn{\omega(E_\alpha)})}\right) + 2^{\ell+1} n\exp\left(-\frac{g}{24(\ell-1)^2}\right).
\end{align*}
\COMMENT{$\ell-1 + n(\ell-1) + n(2^\ell-2) \le n2^{\ell+1}$ since $\ell-1-2n\le 0$ and $\ell-1\le 2^\ell$}
Since $(\ell-1)^{-2}-\ell^{-2}\ge \ell^{-3}$ and $g/24\ell^3\ge \log (2^{\ell+1} n)$ by assumption, we can finally conclude that\COMMENT{if $\tau \ge n$, then $\omega(E_\alpha)=X_n=Y_n$ and $\expn{\omega(E_\alpha)}=X_0=Y_0$}
\begin{align*}
\prob{|\omega(E_\alpha)-  \expn{\omega(E_\alpha)}|> \lambda }
&\leq \prob{|Y_n-Y_0|\geq \lambda}+ \prob{\tau<n}\\
&\leq  2^{\ell}\exp\left(- \frac{\lambda^2}{12\ell ^2M_\ell(\lambda + \expn{\omega(E_\alpha)})}\right) + \exp\left(-\frac{g}{24\ell^2}\right).
\end{align*}
This completes the proof.
\endproof

We are now ready to prove Theorem~\ref{thm:hypermatching2}. The proof proceeds in three steps as outlined in the beginning of this section.

\lateproof{Theorem~\ref{thm:hypermatching2}}

\begin{step}
	Random vertex partition
\end{step}

Let $p:=\Delta^{20Lr\eps  }$.
We will first partition $V(\cH)$ into $p$ subsets $V_1,\dots,V_p$. For each $i\in[p]$, let $\cH_i:=\cH[V_i]$. For an edge $e\in E(\cH)$, let $\tau(e)=i$ if $e\in E(\cH_i)$, and let $\tau(e)=0$ if no such $i$~exists. For a tuple $\cE=(e_1,\dots,e_\ell)\in \binom{E(\cH)}{\ell}$, define the multiset $\tau(\cE):=\Set{\tau(e_1),\dots,\tau(e_\ell)}$.
Let $\cJ_\ell$ be the set of all multisets of size $\ell$ with elements in~$[p]$. For $J\in \cJ_\ell$, let $\supp(J)$ be the underlying set. We further define $\pi(J)$ as the number of functions $f\colon [\ell]\to \supp(J)$ with $\Set{f(1),\dots,f(\ell)}=J$.\COMMENT{This is $=1$ if $J$ has only one $i$, and $\ell!$ if $J$ is crossing. Before: $\pi(J):=\binom{\ell}{a_1,\dots,a_p}$, where $a_1,\dots,a_p\ge 0$ are the multiplicities of the elements of $[p]$ in~$J$.}
For all $\ell\in [L]$ and $J\in \cJ_\ell$, we define $E_J$ as the set of all $\cE\in \binom{E(\cH)}{\ell}$ with $\tau(\cE)=J$.

We claim that there exists a partition $V_1,\dots,V_p$ of $V(\cH)$ such that the following hold:
\begin{enumerate}[label=\rm{(\alph*)}] 
	\item $\Delta(\cH_i)\le (1+\Delta^{-2\eps}) \Delta /p^{r-1}$ for all $i\in[p]$;\label{slicing:max deg}
	\item $\omega(E_J) = (1\pm \Delta^{-2\eps})\omega(E(\cH))\frac{\pi(J)}{p^{r\ell }}$ for all $\ell\in [L]$, $\omega\in \cW_\ell$ and $J\in \cJ_\ell$.\label{slicing:weights}
\end{enumerate}
This can be seen using a probabilistic argument. For every vertex $x\in V(\cH)$ independently, choose an index $i\in[p]$ uniformly at random and assign $x$ to~$V_i$. We now show that \ref{slicing:max deg} and \ref{slicing:weights} hold with high probability, implying that such a partition exists.

For \ref{slicing:max deg}, consider a vertex $x\in V(\cH)$ and $i\in[p]$. 
Let $X$ be the number of edges $e$ containing $x$ for which $e\sm \Set{x}\In V_i$.
For each edge $e$ containing~$x$, we have that $\prob{e\sm \Set{x} \In V_i} = (1/p)^{r-1}$. 
Thus, $\expn{X}=\dg_\cH(x)/p^{r-1}\le \Delta/p^{r-1}$. Note that for any other vertex $y\neq x$, the random label that we choose for $y$ affects $X$ by at most $\dg_\cH(xy)\le \Delta^c(\cH)$. Note that $\sum_{y\in V(\cH)\sm\Set{x}}\dg_\cH(xy) = \dg_\cH(x) (r-1)\le \Delta r$. Thus, using McDiarmid's inequality,\COMMENT{changing the label $x$ doesn't affect $X$ at all} we deduce that
\begin{align*}
\prob{X -\expn{X} \ge \Delta^{1-2\eps}/p^{r-1}} 
&\le 2\exp\left({-\frac{2\Delta^{2-4\eps}}{\Delta^c(\cH) \Delta r p^{2r-2}}}\right) 
\le 2\exp\left(-\Delta^{\delta-45Lr^2\eps}\right)
\leq \exp(-\Delta^\eps).
\end{align*}
With a union bound over all (non-isolated) vertices\COMMENT{There are at most $re(\cH)\le r\exp(\Delta^{\eps^2})$ non-isolated vertices} and $i\in[p]$, we can infer that with high probability~\ref{slicing:max deg} holds.

For \ref{slicing:weights}, consider $\ell\in [L]$, $\omega\in \cW_\ell$ and $J\in \cJ_\ell$. For an edge $e\in E(\cH)$ and $i\in[p]$, we have that $\prob{e\in E(\cH_i)}=p^{-r}$. Thus, for $\cE\in \binom{E(\cH)}{\ell}$, we have $\prob{\tau(\cE)=J}= \pi(J)p^{-r\ell}$ if the edges in $\cE$ are pairwise disjoint, and $\omega(\cE)=0$ otherwise since $\omega$ is clean. Hence, $\expn{\omega(E_J)} = \omega(E(\cH))\frac{\pi(J)}{p^{r\ell}}$. We now establish concentration.
For any vertex~$x$, the random label chosen for $x$ affects $\omega(E_J)$ by at most $\omega(E_x^\ell)$, where $E_x^\ell$ is the set of all $\cE\in \binom{E(\cH)}{\ell}$ for which $x$ is contained in some edge of~$\cE$. Note that 
\begin{align*}
\omega(E_x^\ell)\le \Delta \normV{\omega}{1} \mbox{ for all } x\in V(\cH), \mbox{ and }\sum_{x\in V(\cH)}\omega(E_x^\ell) = r\ell  \omega(E(\cH)).
\end{align*}
Thus, we can use McDiarmid's inequality to conclude that
\begin{align*}
	\prob{\omega(E_J)\neq (1\pm \Delta^{-2\eps})\expn{\omega(E_J)}}  &\le   2\exp\left(-\frac{2\expn{\omega(E_J)}^2}{\Delta \normV{\omega}{1}  r\ell\omega(E(\cH))\Delta^{4\eps}}\right) 
	\le  2\exp\left(-\frac{\omega(E(\cH))}{\normV{\omega}{1} \Delta^{1+45L^2r^2\eps}}\right) \\
	&\le 2\exp\left(-\Delta^{\delta-45L^2r^2\eps}\right) \le \exp(-\Delta^\eps),
\end{align*}
which together with a union bound over all $\ell\in [L]$, $\omega\in \cW_\ell$ and $J\in \cJ_\ell$ proves~\ref{slicing:weights}.

\begin{step}
	Random edge partition
\end{step}

Let $\cH':=\bigcup_{i\in [p]}\cH_i$. 
For each $i\in[p]$, we now partition $\cH_i$ further into $q:=\Delta^{1-20(r-1+1/4L)Lr\eps}$ edge-disjoint subgraphs $\cH_{i,1},\dots,\cH_{i,q}$. Note that
\begin{align}
 	 p^{r-1}q = \Delta^{1-5r\eps}	 \mbox{ and } p^rq \ge \Delta^{1+15Lr\eps}. \label{pq}
\end{align}
We do so (for all $i$ at once) by choosing a function $f\colon E(\cH')\to [q]$ and then let $\cH_{i,j}$ consist of all edges $e\in E(\cH_i)$ with $f(e)=j$, for all $i\in[p],j\in[q]$. 

For $\ell\in [L]$, $J\in \cJ_{\ell}$ and a function $\sigma\colon \supp(J)\to [q]$, let $E_{J,\sigma}$ be the set of all $\cE\in E_J$ for which $\sigma(\tau(e))=f(e)$ for all $e\in \cE$.

We claim that there exists a choice of $f$ such that the following hold:
\begin{enumerate}[label=\rm{(\Alph*)}] 
	\item $\Delta(\cH_{i,j})\le (1+2\Delta^{-2\eps}) \Delta/qp^{r-1}$ for all $i\in [p],j\in[q]$;\label{edge slicing:max deg}
	\item $\Delta^c(\cH_{i,j})\le \Delta^{\eps}$ for all $i\in [p],j\in[q]$;\label{edge slicing:max codeg}
	\item $\omega(E_{J,\sigma})\leq 2\ell!\omega(E(\cH))/q^\ell p^{r\ell}$ for all $\ell\in [L]$, $\omega\in \cW_\ell$, $J\in \cJ_{\ell}$ and $\sigma\colon \supp(J)\to [q]$.\label{edge slicing:weights}
\end{enumerate}

This again can be seen using a probabilistic argument. For each $e\in E(\cH')$ independently, choose $f(e)\in [q]$ uniformly at random.

For \ref{edge slicing:max deg}, fix $i\in [p],j\in[q]$ and a vertex $x\in V(\cH_i)$. Note that $\expn{\dg_{\cH_{i,j}}(x)} = \dg_{\cH_i}(x)/q \le (1+\Delta^{-2\eps})\Delta/qp^{r-1}$ by~\ref{slicing:max deg}. Thus, by Chernoff's bound, we have
\begin{align*}
\prob{\dg_{\cH_{i,j}}(x) -\expn{\dg_{\cH_{i,j}}(x)} \ge \Delta^{1-2\eps}/qp^{r-1}} 
&\le 2\exp\left({-\frac{\Delta^{1-2\eps}}{3q p^{r-1}}}\right) 
\overset{\eqref{pq}}{\leq} \exp\left({-\Delta^{\eps}}\right).
\end{align*}

Similarly, for \ref{edge slicing:max codeg}, fix $i\in [p],j\in[q]$ and two distinct vertices $x,y\in V(\cH_i)$. Note that $\expn{\dg_{\cH_{i,j}}(xy)} = \dg_{\cH_i}(xy)/q \le \Delta^c(\cH)/q\le 1$. Thus, by Chernoff's bound, we have
\begin{align*}
\prob{\dg_{\cH_{i,j}}(xy) \ge \Delta^{\eps}} &\le 2\exp\left({-\Delta^{\eps}}\right).
\end{align*}

To prove \ref{edge slicing:weights}, consider $\ell\in [L]$, $\omega\in \cW_\ell$, $J\in \cJ_{\ell}$ and $\sigma\colon \supp(J)\to [q]$. First note that $\expn{\omega(E_{J,\sigma})} = \omega(E_J)/q^\ell \le \frac{3}{2} \ell!\omega(E(\cH))/q^\ell p^{r\ell}$ by~\ref{slicing:weights}.
We now aim to employ Theorem~\ref{thm: edge slicing} with $E(\cH')$ playing the role of~$V$. 
Let $<$ be an ordering of $E(\cH')$ in which the edges of $\cH_i$ precede those of $\cH_{i'}$ whenever $i<i'$. Write $J=\Set{j_1,\dots,j_\ell}$ such that $j_1\le \dots \le j_\ell$ and define $\alpha:=(\sigma(j_1),\dots,\sigma(j_\ell))\in [q]^\ell$. Hence, for $\cE\in E_J$, we have $\cE\in E_{J,\sigma}$ if and only if $f(e_i)=\sigma(j_i)$ for all $i\in[\ell]$, where $\cE=\Set{e_1,\dots,e_\ell}$ with $e_1<\dots<e_\ell$.\COMMENT{Observe that equal elements in $J$ have equal elements under $\sigma$ so the ordering inside $\cH_i$ is irrelevant.}
Consequently, with notation as in Theorem~\ref{thm: edge slicing}, we have $E_{J,\sigma}=E_J\cap E_\alpha$. 
Thus, $\omega(E_{J,\sigma})=\omega_J(E_\alpha)$, where $\omega_J(\cE):=\omega(\cE)\IND_{\cE\in E_J}$.

We now apply Theorem~\ref{thm: edge slicing} with $E(\cH')$, $\ell$, $\omega_J$, $\frac{1}{2} \ell!\omega(E(\cH))/q^\ell p^{r\ell}$, $\Delta^{2\eps}$ playing the roles of $V,\ell,\omega,\lambda,g$, respectively. For $k\in [\ell]$, we have that (recall that $\omega(E(\cH))\ge \normV{\omega}{k}\Delta^{k+\delta}$ by assumption)
\begin{align*}
\normV{\omega_J}{k}q^{k}g^{k-1} 
\le \normV{\omega}{k} \Delta^{k} 
\le \omega(E(\cH))\Delta^{- \delta}.
\end{align*}\COMMENT{Since $2\eps \le 20(r-1+1/4L)Lr\eps$}
Hence, we infer that (note $\expn{\omega_J(E_\alpha)} +\lambda \le 2\lambda$)
\begin{align*}
\prob{\omega_J(E_\alpha)\ge \expn{\omega_J(E_\alpha)} + \lambda }
&\leq 2^{\ell}\exp\left(- \frac{\lambda}{24 \ell^2 q^{-\ell}\omega(E(\cH))\Delta^{- \delta}}\right) + \exp\left(-\frac{\Delta^{2\eps}}{24\ell^2}\right) \\
&\le 2^{\ell}\exp\left(- \frac{\Delta^{\delta}}{48 \ell p^{r\ell}}\right) + \exp\left(-\frac{\Delta^{2\eps}}{24\ell^2}\right) \le \exp\left({-\Delta^{\eps}}\right).
\end{align*}

A union bound implies that the random choice of $f$ satisfies \ref{edge slicing:max deg}, \ref{edge slicing:max codeg} and \ref{edge slicing:weights} simultaneously with positive probability. From now, fix such a function~$f$.

\begin{step}
	Random matchings
\end{step}

Let $\widetilde{\Delta}:= (1+2\Delta^{- 2\eps})\Delta/qp^{r-1} \ge \Delta^{5r\eps}$ by~\eqref{pq} and $M:=(1+\Delta^{- 2\eps})\widetilde{\Delta}$. Note that
\begin{align}
	p^{r-1}qM = (1\pm 4\Delta^{-2\eps})\Delta.  \label{pqM}  
\end{align}

By~\ref{edge slicing:max deg}, we have $\Delta(\cH_{i,j}) \le \widetilde{\Delta}$.
Moreover, by~\ref{edge slicing:max codeg}, $\Delta^c(\cH_{i,j}) \le \Delta^{\eps}\le \widetilde{\Delta}^{1/5r}$.
Thus, for all $i\in[p]$, $j\in[q]$, we can apply Theorem~\ref{thm:MR} (with $\delta=1/2$, say) to obtain a partition of $E(\cH_{i,j})$ into $M$~matchings.\COMMENT{Theorem~\ref{thm:MR} with $\delta=1/2$ gives $\widetilde{\Delta}+ \widetilde{\Delta}^{1-1/2r}\log^5 \widetilde{\Delta}$ matchings. Want this to be at most $M=\widetilde{\Delta}+\Delta^{-2\eps}\widetilde{\Delta}$ which is fine since $\Delta^{2\eps}\le \widetilde{\Delta}^{2/5r}\le \widetilde{\Delta}^{1/2r}/\log^5\widetilde{\Delta}$.}
This yields a partition of each $E(\cH_i)$ into $q\cdot M$ matchings $\cM_{i,1},\dots,\cM_{i,qM}$.

Now, for each $i\in[p]$ independently, pick an index $s_i\in [qM]$ uniformly at random, and define $$\cM:=\bigcup_{i\in[p]}\cM_{i,s_i}.$$ 
Clearly, $\cM$ is a matching in~$\cH'\In \cH$. Moreover, every edge of $\cH'$ belongs to $\cM$ with probability~$1/qM$. 

Now, consider $\ell\in [L]$ and $\omega\in \cW_\ell$. We first determine the expected value of~$\omega(\cM)$. By linearity, $$\expn{\omega(\cM)} = \sum_{\cE\in \binom{E(\cH)}{\ell}} \omega(\cE)\prob{\cE\In \cM}.$$
We analyse this sum according to the different types of~$\cE$.
For $k\in[\ell]$, let $\cJ_{\ell,k}$ be the set of all $J\in \cJ_\ell$ with $|\supp(J)|=k$. 
Consider $\cE\in  \binom{E(\cH)}{\ell}$ and let $J:=\tau(\cE)$. Note that if $0\in J$, then some edge in $\cE$ does not belong to $\cH'$ and hence $\prob{\cE\In \cM}=0$. Hence, we can assume that $J\in \cJ_\ell$. If $J\in\cJ_{\ell,\ell}$, then the edges in $\cE$ belong to $\cM$ independently with probability~$1/qM$, and hence $\prob{\cE\In \cM}=(qM)^{-\ell}$. Now, suppose $J\in\cJ_{\ell,k}$ for some $k\in[\ell-1]$. By the definition of~$\cM$, if $e,e'\in \cE$ with $\tau(e)=\tau(e')$, then $\prob{\cE\In \cM}=0$ if $e\in E(\cH_{\tau(e),j})$ and $e'\in E(\cH_{\tau(e),j'})$ for distinct $j,j'$. Hence, we can further assume that $\cE\in E_{J,\sigma}$ for some $\sigma\colon \supp(J)\to [q]$. We then have $\prob{\cE\In \cM}\in \Set{0,(qM)^{-k}}$.\COMMENT{$0$ if the edges in the same edge slice do not end up in the same matching from Molloy--Reed.} Altogether, we deduce that
\begin{align*}
\expn{\omega(\cM)} &= \sum_{J\in \cJ_{\ell,\ell}} \omega(E_J)(qM)^{-\ell} \pm \sum_{k=1}^{\ell-1} \sum_{J\in \cJ_{\ell,k},\sigma\colon \supp(J)\to [q]} \omega(E_{J,\sigma})(qM)^{-k}.
\end{align*}
We will show that the first sum is the dominant term. Clearly, $|\cJ_{\ell,\ell}|=\binom{p}{\ell}$. Thus, using \ref{slicing:weights}, we infer that
\begin{align*}
\sum_{J\in \cJ_{\ell,\ell}} \omega(E_J)(qM)^{-\ell} 
= \binom{p}{\ell}\cdot (1\pm \Delta^{-2\eps})\frac{\ell!\omega(E(\cH))}{p^{r\ell}}\cdot \frac{1}{(qM)^{\ell}} \overset{\eqref{pqM}}{=} (1\pm \Delta^{-3\eps/2}) \omega(E(\cH)) /\Delta^\ell.
\end{align*}\COMMENT{$\binom{p}{\ell}\ell!=p^\ell(1\pm \O_\ell(1/p))=p^\ell(1\pm \Delta^{-19\eps})$}
For $k\in[\ell-1]$, employing \ref{edge slicing:weights} and $|\cJ_{\ell,k}|=\binom{p}{k}\binom{\ell-1}{k-1}$,\COMMENT{$\binom{p}{k}$ ways to choose the $k$ elements that appear in the multiset. Order them. Then have $\binom{\ell-1}{k-1}$ ways to produce a multiset by putting $k-1$ `separators' in the gaps between $\ell$ consecutive objects.} we deduce that
\begin{align*}
\sum_{J\in \cJ_{\ell,k},\sigma\colon \supp(J)\to [q]} \omega(E_{J,\sigma})(qM)^{-k} \le p^k 2^\ell q^k \frac{2\ell!\omega(E(\cH))}{q^\ell p^{r\ell}} \cdot \frac{1}{(qM)^{k}}   \le \frac{\omega(E(\cH))}{\Delta^{\ell+14\eps}},
\end{align*}
\COMMENT{$q^k$ number of $\sigma$}where in the last inequality we used that $\frac{p^k q^k}{q^\ell p^{r\ell}(qM)^{k}} = \frac{1}{(p^rq)^{\ell-k}(p^{r-1}qM)^k}$ together with $(p^rq)^{\ell-k}\ge \Delta^{\ell-k+15\eps}$ by \eqref{pq} and $(p^{r-1}qM)^k\ge \frac{1}{2}\Delta^k$ by~\eqref{pqM}.\COMMENT{We have $(p^rq)^{\ell-k}\ge \Delta^{(\ell-k)(1+10Lr\eps)}\ge \Delta^{\ell-k+15\eps}$ since $\ell-k\ge 1$. Absorb constants by reducing $15$ to $14$}
Putting everything together, we obtain that $$\expn{\omega(\cM)} = (1\pm 2\Delta^{-3\eps/2}) \omega(E(\cH)) /\Delta^\ell.$$

Finally, we need to bound the effect of each random variable~$s_i$. Note that each outcome of the variables $s_1,\dots,s_p$ induces a function $\sigma\colon [p]\to [q]$, where $\sigma(i)$ is the unique $j\in[q]$ for which $\cM_{i,s_i}$ was one of the matchings coming from $E(\cH_{i,j})$,\COMMENT{before: $\cM_{i,s_i}\In E(\cH_{i,j})$, but $\cM_{i,s_i}$ could be empty and then $j$ would not be unique} and each tuple $\cE\In \cM$ satisfies $\cE\in E_{J,\sigma|_{\supp(J)}}$, where $J=\tau(\cE)\in \cJ_\ell$. Since changing the value of $s_i$ only affects those $\cE$ with $i\in \tau(\cE)$, we have that the effect of $s_i$ on $\omega(\cM)$ is at most
\begin{align*}
\max_{\sigma\colon [p]\to [q]} \sum_{J\in \cJ_\ell\colon i\in J} \omega(E_{J,\sigma|_{\supp(J)}})  \overset{\ref{edge slicing:weights}}{\leq} p^{\ell-1}\frac{2\ell!\omega(E(\cH))}{q^\ell p^{r\ell}} \overset{\eqref{pq}}{=}\frac{2\ell!\omega(E(\cH))}{p\Delta^{(1-5r\eps)\ell}} \le \frac{\omega(E(\cH))}{\Delta^{\ell+14Lr\eps}}.
\end{align*}
\COMMENT{number of $J$ with $i\in J$ is at most $p^{\ell-1}$. last exponent of $\Delta$: $20Lr\eps$ from $p$ minus $5r\ell \eps$ leaves at least $15Lr\eps$. reduce by one to absorb constants}

Thus, using McDiarmid's inequality, we deduce that
\begin{align*}
\prob{\omega(\cM)\neq (1\pm \Delta^{-2 \eps}) \expn{\omega(\cM)}} 
&\le 2\exp\left({-\frac{2\Delta^{-4 \eps} \expn{\omega(\cM)}^2}{p\cdot (\omega(E(\cH))/\Delta^{\ell+14Lr\eps})^2}}\right) \\
&\le 2\exp\left({-\frac{ \Delta^{28Lr\eps -4 \eps}}{p}}\right)
\leq  
\exp\left(-\Delta^{\eps} \right).
\end{align*}
A union bound over all $\ell\in [L]$ and $\omega\in \cW_\ell$ completes the proof.
\endproof

\section{Applications}\label{sec:apps}

In this section we provide a small exposition of applications of Theorem~\ref{thm:hypermatching2}. In Section~\ref{subsec:pseudo Steiner},
we deduce the existence of approximate Steiner systems that behave `randomly', e.g.~with respect to subgraph statistics.
Then, we briefly explain how we apply it in two forthcoming papers~\cite{EGJ:19b, EJ:19}
on rainbow embeddings and approximate decompositions.

\subsection{Pseudorandom Steiner systems} \label{subsec:pseudo Steiner}
Recall that an $(n,k,t)$-Steiner system is a set $\cS$ of $k$-subsets of some $n$-set $V$ such that every $t$-subset of $V$ is contained in exactly one $k$-set in~$\cS$. We now view such $\cS$ as a $k$-graph. Note that any subgraph of $\cS$ has the following property: any two of its edges intersect in less than $t$~vertices; we will simply say that such graphs are \defn{$t$-avoiding}. For $t=2$, such hypergraphs are often called `linear' or `simple'. Now, for a fixed $t$-avoiding $k$-graph $F$, we may ask how many (labelled) copies of $F$ exist in~$\cS$. Since $|\cS|=\binom{n}{t}/\binom{k}{t}$, the edge density of $\cS$ is (for large~$n$) approximately $p:=(k-t)!n^{-k+t}$. In a random $k$-graph with this density, we would expect $p^{e(F)}n^{v(F)}$ labelled copies of~$F$. Of course this makes only sense when $(-k+t)e(F)+v(F)>0$, or equivalently, when the average degree of $F$ is less than~$k/(k-t)$. Moreover, in order to be able to obtain precise counts for~$F$, one needs this condition for all non-empty subgraphs of~$F$. We thus define the \defn{maximum average degree of $F$}, denoted $\mad(F)$, as the maximum of $ke(F')/v(F')$ over all non-empty subgraphs $F'$ of~$F$. For two $k$-graphs $F,G$, let $\inj(F,G)$ be the number of labelled copies of $F$ in $G$, that is, the number of injections $f\colon V(F)\to V(G)$ for which $f(e)\in E(G)$ for all $e\in E(F)$.

\COMMENT{There are various problems concerning the count of subconfigurations in Steiner systems, and more general designs. Erdos conjecture on sparse Steiner triple systems, where some counts should be $0$, for example Pasch configuration. Others want to maximize the number of such configurations. Stinson and Wei}
As one application of Theorem~\ref{thm:hypermatching2}, we show that there exist approximate Steiner systems whose subgraph statistics resemble the random model.

\begin{theorem}
Suppose $1/n\ll \eps \ll 1/k,1/v$ and $t\in \Set{2,\dots,k-1}$. Let $\cF$ be the family of all $t$-avoiding $k$-graphs $F$ with $v(F)\le v$ and $\mad(F)<k/(k-t)$, and let $p:=(k-t)!n^{-k+t}$.
There exists a partial $(n,k,t)$-Steiner system $\cS$ with $|\cS|\ge (1-n^{-\eps})\binom{n}{t}/\binom{k}{t}$ such that
\begin{align*}
	\inj(F,\cS) = (1\pm n^{-\eps})p^{e(F)} n^{v(F)} \mbox{ for all }F\in \cF.
\end{align*}
\end{theorem}

\proof
Choose a new constant $\delta>0$ such that $1/n\ll \eps \ll \delta \ll 1/k,1/v$.\COMMENT{$\delta\approx 1/(k-t)$ should suffice}

For $e\in \binom{[n]}{k}$, let $\pi(e):=\binom{e}{t}$, and for a $k$-graph $G$ on~$[n]$, let $\pi(G)$ be the $\binom{k}{t}$-uniform hypergraph on $V:=\binom{[n]}{t}$ with edge set $\set{\pi(e)}{e\in E(G)}$.  
Let $\cH:=\pi(K_n^k)$, which is the hypergraph already defined in Section~\ref{subsec:pippenger}. Note that $\pi$ is a bijection between $k$-graphs on $[n]$ and spanning subgraphs of~$\cH$.

Recall that $\cH$ is $\binom{n-t}{k-t}$-regular. Moreover, $\Delta^c(\cH)\le n^{k-t-1} \le \Delta(\cH)^{1-\delta}$. Crucially, matchings in $\cH$ correspond to $t$-avoiding subgraphs of $K_n^k$, and thus to partial $(n,k,t)$-Steiner systems.

Now, fix $F\in \cF$ and let $\ell:=e(F)$. Define the $\ell$-tuple weight function $\omega_F$ on $E(\cH)$ as follows: for an $\ell$-set $\cE=\Set{\pi(e_1),\dots,\pi(e_\ell)}$ of edges of~$\cH$, let $\omega_F(\cE)$ be the number of injections $f\colon V(F)\to [n]$ for which $\set{f(e)}{e\in E(F)}=\Set{e_1,\dots,e_\ell}$. Hence, for $G\In K_n^k$, we have  $\inj(F,G)=\omega_F(\pi(G))$. In particular, $\omega_F(E(\cH))=(1\pm n^{-0.9})n^{v(F)}$. Note also that $\omega_F$ is clean since $F$ is $t$-avoiding.

Fix $\ell'\in[\ell]$ and a set of $\ell'$ edges $\pi(e_1),\dots,\pi(e_{\ell'})$ in~$\cH$. Let $v':=|e_1\cup \dots \cup e_{\ell'}|$. The number of injections $f \colon V(F)\to [n]$ for which $\Set{e_1,\dots,e_{\ell'}}\In \set{f(e)}{e\in E(F)}$ is at most $v(F)!n^{v(F)-v'}$.
Since $\mad(F)<k/(k-t)$, we have $v'>\ell'(k-t)$, implying
\begin{align*}
	\normV{\omega_F}{\ell'}\Delta^{\ell'+\delta} \le v! n^{v(F)-v'}\cdot n^{(k-t)(\ell'+\delta)} \le n^{v(F)}/2 \le \omega_F(E(\cH)).
\end{align*}

Thus, we can apply\COMMENT{only constantly many weight functions, $e(\cH)=\binom{n}{k}\le \exp(n^{\eps^2})$} Theorem~\ref{thm:hypermatching2} with $\Delta:=1/p=n^{k-t}/(k-t)!\ge \Delta(\cH)$ to obtain a matching~$\cM$ such that $\omega_F(\cM)=(1\pm \Delta^{-\eps}) \omega_F(E(\cH))/\Delta^{e(F)}$ for all $F\in \cF$. Let $\cS:=\pi^{-1}(\cM)$. Note that $\cS$ is a partial $(n,k,t)$-Steiner system. Moreover, for any $F\in \cF$, we have $$\inj(F,\cS)=\omega_F(\cM)=(1\pm \Delta^{-2\eps})(1\pm n^{-0.9})n^{v(F)}p^{e(F)}=(1\pm n^{-\eps})p^{e(F)} n^{v(F)},$$ as desired.

Finally, note that the $k$-graph $F_0$ consisting of only one edge is trivially $t$-avoiding and $\mad(F_0)=1$. Thus, by the above, $\inj(F_0,\cS)\ge (1-n^{-\eps})p n^k$. We conclude that
$$|\cS|=\inj(F_0,\cS)/k! \ge (1-n^{-\eps})p n^k/k! \ge (1-n^{-\eps})\binom{n}{t}/\binom{k}{t},$$
completing the proof.
\endproof

One could also easily ensure that the residual $t$-graph of uncovered $t$-sets is quasirandom. This would allow for an application of a clique decomposition result to complete $\cS$ to a Steiner system. The lower bound on the number of $F$-copies would then still hold. However, such a completion step, even if only applied to $o(n^t)$ $t$-sets, could drastically increase the number of $F$-copies. For simplicity, we thus omitted such a completion entirely. It is needless to say that variations of this theorem can be obtained in the same way, for instance asking for the number of `rooted' copies.

\subsection{Rainbow problems}\label{subsec:rainbow}

In~\cite{EGJ:19b}, we consider subgraph embeddings in edge-coloured graphs with the additional requirement that the embedded subgraph is `rainbow', meaning that any two edges in the subgraph have distinct colours. Such rainbow embeddings have applications to various other problems.
For instance, Montgomery, Pokrovskiy and Sudakov~\cite{MPS:18} recently used rainbow embeddings to approximately solve Ringel's conjecture from 1963 (which states that any tree with $n$ edges decomposes~$K_{2n+1}$).
We consider the classic setting of the blow-up lemma due to Koml\'os, Sark\"ozy and Szemer\'edi~\cite{KSS:97}. Given a multipartite graph $G$, where the bipartite graphs between two parts are `quasirandom', and a bounded degree graph $H$ with a fitting vertex partition, $H$ can be embedded as a spanning subgraph of~$G$. We show that this is still true when $G$ is edge-coloured and we want to find a rainbow copy of~$H$, assuming certain boundedness conditions on the edge-colouring which can be seen to be almost optimal.

To achieve this, we employ Theorem~\ref{thm:hypermatching simple} as a crucial tool.
In the following, we briefly explain how we apply Theorem~\ref{thm:hypermatching simple} and exploit the weight functions in our proof.
To this end, we consider the following toy example.
Suppose $G$ is the complete bipartite graph with bipartition $(U,V)$ and $|U|=|V|=n$.
Suppose further that $c\colon E(G) \to C$ is a proper edge-colouring of~$G$.
Our aim is to find an almost perfect rainbow matching.
When the colouring is optimal, then finding such a matching of size $n-1$ is equivalent to the famous Ryser--Brualdi--Stein conjecture on almost transversals in Latin squares.

In order to apply our theorem, we formulate the problem as a hypergraph matching problem. Let $\cH$ be the hypergraph with vertex set $U\cup V\cup C$ and edge set $\{\{u,v,c(uv)\}\colon uv\in E(G)\}$.
The key property of $\cH$ is the following bijection between the set of all rainbow matchings in~$G$ and the set of all matchings in $\cH$ --
we simply assign a rainbow matching $M$ in $G$ to the matching $\cM:=\{\{u,v,c(uv)\}\colon uv\in M\}$ in~$\cH$.
Clearly, $\Delta(\cH)=n$ and $\Delta^c(\cH)=1$. The existence of an almost perfect rainbow matching in~$G$ follows now from known results. For instance, Theorem~\ref{thm:MR} yields a decomposition of~$E(\cH)$ into $(1+o(1))n$ hypergraph matchings, and as $e(\cH)=n^2$, there must be a hypergraph matching~$\cM$ of size $(1-o(1))n$ in $\cH$ and in turn a rainbow matching $M$ in $G$ of this size.

In the proof of our rainbow blow-up lemma, we also seek almost perfect matchings in bipartite graphs. However, we need much more control over these matchings, which we achieve using our new Theorem~\ref{thm:hypermatching simple}. Our proof proceeds in several rounds where in each round, we embed essentially all vertices that need to be embedded into a particular cluster of our multipartite graph. Each such embedding step is modelled as finding a rainbow matching $M$ in an auxiliary bipartite `candidacy graph'. Although these candidacy graphs are more complicated and have more complex colour constraints than our toy example above, they can still be handled using hypergraph matchings in a similar way.
However, in order to perform the embedding rounds repeatedly, we need to ensure that certain quasirandomness properties are preserved throughout the procedure, which depend on the previous embeddings.
In our toy example this would mean, for instance, that for some specified sets $U'\In U,V'\In V$, we need $|E(G[U',V'])\cap M|\approx |U'||V'|/n$,
and more generally that for sets $E'\In E(G)$, we need $|E'\cap M|\approx |E'|/n$.
This can be ensured by utilizing weight functions as in Theorem~\ref{thm:hypermatching simple} by defining $\omega_{E'}(e\cup \Set{c(e)})=\IND_{e\in E'}$ for all $e\in E(G)$.
It is very important here that Theorem~\ref{thm:hypermatching simple} applies to hypergraphs which are not necessarily almost regular, since the colouring can be arbitrary. For the same reason, it is useful that $v(\cH)$ plays no role in the parametrization of the theorem.

\subsection{Decompositions}
Kim, K\"uhn, Osthus and Tyomkyn~\cite{KKOT:19} proved that in the setting of the original blow-up lemma as described above,
the quasirandom multipartite graph does not only contain any single graph of bounded degree with the same multipartite structure, but can even be almost decomposed into any collection of such bounded degree graphs. This result has already found fruitful applications~\cite{CKKO:19,JKKO:ta}.
The first and third author give an alternative and in particular much shorter proof for this decomposition result in~\cite{EJ:19}.

The overall strategy is to use Theorem~\ref{thm:hypermatching2} in a similar way as for the rainbow embeddings described in Section~\ref{subsec:rainbow}.
To this end, the decomposition problem is transformed into a rainbow embedding problem as follows.
Suppose $G$ is a graph and $\cH$ is a collection of graphs on at most~$|V(G)|$ vertices.
Define a new graph $\cG$ by taking $|\cH|$ disjoint copies $(G_H)_{H\in \cH}$ of $G$ and colour every copy of a particular edge of~$G$ with a unique colour.
Hence, a collection of edge-disjoint copies of the graphs in~$\cH$ in~$G$ is equivalent to a rainbow embedding of the disjoint union of the graphs in $\cH$ into $\cG$ where each $H\in \cH$ is embedded into~$G_H$.

\bibliographystyle{../amsplain_v2.0customized}
\bibliography{../ReferencesLocal}

\providecommand{\bysame}{\leavevmode\hbox to3em{\hrulefill}\thinspace}
\providecommand{\MR}{\relax\ifhmode\unskip\space\fi MR }
\providecommand{\MRhref}[2]{%
  \href{http://www.ams.org/mathscinet-getitem?mr=#1}{#2}
}
\providecommand{\href}[2]{#2}
\begin{thebibliography}{10}

\bibitem{AGS:09}
R.~Aharoni, A.~Georgakopoulos, and P.~Spr\"{u}ssel, \emph{Perfect matchings in
  {$r$}-partite {$r$}-graphs}, European J. Combin.~\textbf{30} (2009), 39--42.

\bibitem{AKS:81}
M.~Ajtai, J.~Koml\'{o}s, and E.~Szemer\'{e}di, \emph{A dense infinite {S}idon
  sequence}, European J. Combin.~\textbf{2} (1981), 1--11.

\bibitem{AY:05}
N.~Alon and R.~Yuster, \emph{On a hypergraph matching problem}, Graphs
  Combin.~\textbf{21} (2005), 377--384.

\bibitem{CKKO:19}
P.~Condon, J.~Kim, D.~K\"uhn, and D.~Osthus, \emph{A bandwidth theorem for
  approximate decompositions}, Proc. Lond. Math. Soc.~\textbf{118} (2019),
  1393--1449.

\bibitem{CFMR:96}
C.~Cooper, A.~Frieze, M.~Molloy, and B.~Reed, \emph{Perfect matchings in random
  {$r$}-regular, {$s$}-uniform hypergraphs}, Combin. Probab. Comput.~\textbf{5}
  (1996), 1--14.

\bibitem{EGJ:19b}
S.~Ehard, S.~Glock, and F.~Joos, \emph{A rainbow blow-up lemma for almost
  optimally bounded edge-colourings}, preprint (2019).

\bibitem{EJ:19}
S.~Ehard and F.~Joos, \emph{A short proof of the blow-up lemma for approximate
  decompositions}, preprint (2019).

\bibitem{FR:85}
P.~Frankl and V.~R\"{o}dl, \emph{Near perfect coverings in graphs and
  hypergraphs}, European J. Combin.~\textbf{6} (1985), 317--326.

\bibitem{freedman:75}
D.~A.~Freedman, \emph{On tail probabilities for martingales}, Ann.
  Probab.~\textbf{3} (1975), 100--118.

\bibitem{FK:12}
A.~Frieze and M.~Krivelevich, \emph{Packing {H}amilton cycles in random and
  pseudo-random hypergraphs}, Random Structures Algorithms~\textbf{41} (2012),
  1--22.

\bibitem{furedi:88}
Z.~F\"{u}redi, \emph{Matchings and covers in hypergraphs}, Graphs
  Combin.~\textbf{4} (1988), 115--206.

\bibitem{GKLO:16}
S.~Glock, D.~K\"uhn, A.~Lo, and D.~Osthus, \emph{The existence of designs via
  iterative absorption}, arXiv:1611.06827 (2016).

\bibitem{HPS:09}
H.~H\`an, Y.~Person, and M.~Schacht, \emph{On perfect matchings in uniform
  hypergraphs with large minimum vertex degree}, SIAM J. Discrete
  Math.~\textbf{23} (2009), 732--748.

\bibitem{JKV:08}
A.~Johansson, J.~Kahn, and V.~Vu, \emph{Factors in random graphs}, Random
  Structures Algorithms~\textbf{33} (2008), 1--28.

\bibitem{JKKO:ta}
F.~Joos, J.~Kim, D.~K\"uhn, and D.~Osthus, \emph{Optimal packings of bounded
  degree trees}, J. Eur. Math. Soc. (to appear).

\bibitem{kahn:96}
J.~Kahn, \emph{Asymptotically good list-colorings}, J. Combin. Theory Ser.
  A~\textbf{73} (1996), 1--59.

\bibitem{keevash:14}
P.~Keevash, \emph{The existence of designs}, arXiv:1401.3665 (2014).

\bibitem{keevash:18c}
\bysame, \emph{Hypergraph matchings and designs}, Proc. Int. Cong. of
  Math.~\textbf{3} (2018), 3099--3122.

\bibitem{KM:15}
P.~Keevash and R.~Mycroft, \emph{A geometric theory for hypergraph matching},
  Mem. Amer. Math. Soc.~\textbf{233} (2015), vi+95.

\bibitem{KKOT:19}
J.~Kim, D.~K\"uhn, D.~Osthus, and M.~Tyomkyn, \emph{A blow-up lemma for
  approximate decompositions}, Trans. Amer. Math. Soc.~\textbf{371} (2019),
  4655--4742.

\bibitem{KPS:82}
J.~Koml\'{o}s, J.~Pintz, and E.~Szemer\'{e}di, \emph{A lower bound for
  {H}eilbronn's problem}, J. London Math. Soc. (2)~\textbf{25} (1982), 13--24.

\bibitem{KSS:97}
J.~Koml\'os, G.~N.~S\'ark\"ozy, and E.~Szemer\'edi, \emph{Blow-up lemma},
  Combinatorica~\textbf{17} (1997), 109--123.

\bibitem{mcdiarmid:89}
C.~McDiarmid, \emph{On the method of bounded differences}, Surveys in
  combinatorics, 1989 ({N}orwich, 1989), London Math. Soc. Lecture Note Ser.
  141, Cambridge Univ. Press, 1989, pp.~148--188.

\bibitem{MR:00}
M.~Molloy and B.~Reed, \emph{Near-optimal list colorings}, Random Structures
  Algorithms~\textbf{17} (2000), 376--402.

\bibitem{MPS:18}
R.~Montgomery, A.~Pokrovskiy, and B.~Sudakov, \emph{Embedding rainbow trees
  with applications to graph labelling and decomposition}, J. Eur. Math. Soc.
  (to appear).

\bibitem{PS:89}
N.~Pippenger and J.~Spencer, \emph{Asymptotic behaviour of the chromatic index
  for hypergraphs}, J. Combin. Theory Ser.~A~\textbf{51} (1989), 24--42.

\bibitem{rodl:85}
V.~R\"odl, \emph{On a packing and covering problem}, European J.
  Combin.~\textbf{6} (1985), 69--78.

\bibitem{RR:10}
V.~R\"odl and A.~Ruci\'nski, \emph{Dirac-type questions for hypergraphs---a
  survey (or more problems for {E}ndre to solve)}, In:~An irregular mind
  (I.~B\'ar\'any and J.~Solymosi, eds.), Bolyai Soc. Math. Stud.~21, J\'anos
  Bolyai Math. Soc., 2010, pp.~561--590.

\bibitem{RRS:09}
V.~R\"{o}dl, A.~Ruci\'{n}ski, and E.~Szemer\'{e}di, \emph{Perfect matchings in
  large uniform hypergraphs with large minimum collective degree}, J. Combin.
  Theory Ser. A~\textbf{116} (2009), 613--636.

\end{thebibliography}

\bigskip

\end{document}